%% file: Incomplete_PCM_theory_applications.tex
\documentclass[graybox]{svmult}

\usepackage[T1]{fontenc}

\usepackage{type1cm}        % activate if the above 3 fonts are
\usepackage{makeidx}         % allows index generation
\usepackage{graphicx}        % standard LaTeX graphics tool
\usepackage{multicol}        % used for the two-column index
\usepackage[bottom]{footmisc}% places footnotes at page bottom

\usepackage{newtxtext}       % 
\usepackage[varvw]{newtxmath}       % selects Times Roman as basic font

\makeindex             % used for the subject index
\usepackage[dvipsnames,svgnames,table]{xcolor}                
\usepackage{hyperref}
\usepackage{array}
\newcolumntype{K}[1]{>{\centering\arraybackslash$}p{#1}<{$}}

\usepackage{pgfplots}
\pgfplotsset{every tick label/.append style={font=\footnotesize}}
\pgfplotsset{compat=1.18}
\usepackage{tikz}
\usetikzlibrary{arrows, calc, matrix, patterns, positioning, trees}

\usepackage{multirow}
\usepackage{tabularx}
\usepackage{threeparttable}
\usepackage{booktabs}
\newcolumntype{R}{>{\raggedleft\arraybackslash}X}
\newcolumntype{L}{>{\raggedright\arraybackslash}X}
\newcolumntype{C}{>{\centering\arraybackslash}X}

\begin{document}

\title*{Incomplete pairwise comparison matrices and their applications}
\author{L\'aszl\'o Csat\'o\orcidID{0000-0001-8705-5036} and\\ S\'andor Boz\'oki\orcidID{0000-0003-4170-4613}}
%\author{\href{https://sites.google.com/view/laszlocsato}{L\'aszl\'o Csat\'o} \\ S\'andor Boz\'oki}
\institute{L\'aszl\'o Csat\'o \at HUN-REN Institute for Computer Science and Control (SZTAKI), Hungary, 1111 Budapest, Kende street 13--17., \email{laszlo.csato@sztaki.hun-ren.hu} \at Corvinus University of Budapest (BCE), Hungary, 1093 Budapest, F{\H o}v\'am square 8.
\and S\'andor Boz\'oki \at HUN-REN Institute for Computer Science and Control (SZTAKI), Hungary, 1111 Budapest, Kende street 13--17., \email{bozoki.sandor@sztaki.hun-ren.hu} \at Corvinus University of Budapest (BCE), Hungary, 1093 Budapest, F{\H o}v\'am square 8.}
\maketitle

%\abstract*{Incomplete pairwise comparison matrices are increasingly used to save resources and reduce cognitive load by collecting only a subset of all possible pairwise comparisons. We present their graph representation and some completion methods, including the incomplete eigenvector and incomplete logarithmic least squares methods, as well as a recently suggested lexicographical minimisation of triad inconsistencies. The issue of ordinal violations is discussed for matrices generated by directed acyclic graphs and the best-worst method. We also show a reasonable approach to generalise the inconsistency threshold based on the dominant eigenvalue to the incomplete case, and state some recent results on the optimal order of obtaining pairwise comparisons. The benefits of using incomplete pairwise comparisons are highlighted by several applications.}

\abstract{Incomplete pairwise comparison matrices are increasingly employed to save resources and reduce cognitive load by collecting only a subset of all possible pairwise comparisons. We present their graph representation and some completion algorithms, including the incomplete eigenvector and incomplete logarithmic least squares methods, as well as a lexicographical minimisation of triad inconsistencies. The issue of ordinal violations is discussed for matrices generated by directed acyclic graphs and the best--worst method. We also show a reasonable approach to generalise the inconsistency threshold based on the dominant eigenvalue to the incomplete case, and state recent results on the optimal order of obtaining pairwise comparisons. The benefits of using incomplete pairwise comparisons are highlighted by several applications.}


% MSC class
% Management decision making, including multiple objectives
% 91B06 Decision theory

\section{Introduction} \label{Sec1}

Pairwise comparison matrices are extensively used in multi-criteria decision-making methods. However, obtaining the $n(n-1)/2$ possible pairwise comparisons of $n$ alternatives might be time-consuming, or even infeasible when they come from an external dataset. Collecting all comparisons can also be inefficient, as a sufficient proportion of them may provide a good approximation of the preferences calculated from all possible comparisons \cite{CarmoneKaraZanakis1997}. For example, the best--worst method requires the identification of the best (most desirable, most important) and the worst (least desirable, least important) alternatives in the first step, and the other alternatives are compared only to these particular alternatives \cite{Rezaei2015}.

Incomplete pairwise comparison matrices, which contain some missing/unknown entries, increasingly appear in both the theoretical literature and practical applications. The current chapter gives a concise overview of this field. The mathematical background is presented in Section~\ref{Sec2}. Section~\ref{Sec3} discusses completion and weighting methods, focusing on the extension of the well-known eigenvector and logarithmic least squares (row geometric mean) methods, originally suggested for complete pairwise comparison matrices, to this general setting. Section~\ref{Sec4} deals with the issue of ordinal violations that may emerge in ranking from incomplete pairwise comparison matrices induced by directed acyclic graphs, or by the best--worst method. Inconsistency thresholds and optimal filling patterns for incomplete pairwise comparison matrices are surveyed in Sections~\ref{Sec5} and \ref{Sec6}, respectively. Section~\ref{Sec7} summarises some applications of incomplete pairwise comparisons, and outlines real-world problems that can be addressed by this methodology due to the structure of the data. Finally, Section~\ref{Sec8} concludes.

\section{Preliminaries} \label{Sec2}

Given a non-empty set of $n$ alternatives, an $n \times n$ positive matrix $\mathbf{A} = \left[ a_{ij} \right]$ is called a \emph{pairwise comparison matrix} if $a_{ji} = 1 / a_{ij}$ for all $1 \leq i,j \leq n$. The entry $a_{ij}$ represents the extent to which alternative $i$ is preferred to alternative $j$.

A pairwise comparison matrix is said to be \emph{consistent} if $a_{ik} = a_{ij} a_{jk}$ holds for all $1 \leq i,j,k \leq n$. Otherwise, it is called \emph{inconsistent}. The level of inconsistency is measured by inconsistency indices; see \cite{Brunelli2018} for an overview of them.

Let $\mathbf{A} = \left[ a_{ij} \right]$ be a pairwise comparison matrix and $\lambda_{\max}(\mathbf{A})$ be its dominant eigenvalue. The inconsistency index suggested by \cite{Saaty1977} is
\[
\mathit{CI} (\mathbf{A}) = \frac{\lambda_{\max}(\mathbf{A})-n}{n-1}.
\]
Denote the average inconsistency index $\mathit{CI}$ of a large number of randomly generated pairwise comparison matrices by $\mathit{RI}_n$, with all entries in the upper triangle chosen uniformly and randomly from the set of integers between 1 and 9 and their reciprocals.
The inconsistency ratio of a matrix $\mathbf{A}$ is $\mathit{CR} (\mathbf{A}) = \mathit{CI} (\mathbf{A}) / \mathit{RI}_n$. According to \cite{Saaty1977}, matrices with $\mathit{CR} < 0.1$ can be accepted; otherwise, the pairwise comparisons need to be revised.

\emph{Incomplete} pairwise comparisons allow for some missing entries, indicated by $\ast$ in the following. Formally, matrix $\mathbf{A} = \left[ a_{ij} \right]$ is an incomplete pairwise comparison matrix if, for all $1 \leq i,j \leq n$,
(i) $a_{ij} \neq \ast$ implies $a_{ji} = 1 / a_{ij}$;
(ii) $a_{ij} = \ast$ implies $a_{ji} = \ast$; and
(iii) $a_{ii} = 1$ for all $1 \leq i \leq n$.

The graph representation of known comparisons is a convenient tool to analyse incomplete pairwise comparison matrices \cite{SzadoczkiBozokiTekile2022}.
Let $\mathbf{A} = \left[ a_{ij} \right]$ be an incomplete pairwise comparison matrix. Its associated graph is the undirected graph $G = (V,E)$, where the set of vertices $V$ represents the set of alternatives ($\lvert V \rvert = n$), and the edge set $E$ contains the known comparisons: $(i,j) \in E \iff a_{ij} \neq \ast$ for all $1 \leq i < j \leq n$.

\begin{example} \label{Examp1}
Consider the following incomplete pairwise comparison matrix with two missing entries above the diagonal:
\[
\mathbf{A} = \left[
\begin{array}{K{3em} K{3em} K{3em} K{3em}}
    1     	& a_{12}  	& \ast   	& a_{14} \\
    a_{21}	& 1       	& a_{23}	& \ast   \\
   \ast		& a_{32}	& 1      	& a_{34} \\
    a_{41}	& \ast  	& a_{43}	& 1 \\
\end{array}
\right].
\]

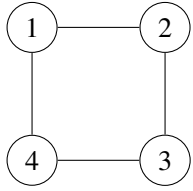
\begin{figure}[t!]
\centering
\begin{tikzpicture}[scale=1.25, auto=center, transform shape, >=triangle 45]
\tikzstyle{every node}=[draw,shape=circle];
  \node (n1) at (135:1) {$1$};
  \node (n2) at (45:1)  {$2$};
  \node (n3) at (315:1) {$3$};
  \node (n4) at (225:1) {$4$};

  \foreach \from/\to in {n1/n2,n1/n4,n2/n3,n3/n4}
    \draw (\from) -- (\to);
\end{tikzpicture}

\caption{The graph representation of the incomplete pairwise comparison matrix in Example~\ref{Examp1}}
\label{Fig1}
\end{figure}

Figure~\ref{Fig1} shows the graph associated with the incomplete matrix $\mathbf{A}$.
\end{example}

The number of pairwise comparisons can be reduced by using reference alternatives, too.
The best--worst method \cite{Rezaei2015} starts by identifying the best and the worst alternatives. In the second step, all other alternatives are compared to the best and the worst, respectively. This approach requires only $2n-3$ comparisons for $n$ alternatives, in contrast to a complete pairwise comparison matrix with $n(n-1)/2$ comparisons. The resulting matrix is called the \emph{best--worst method matrix}: it is an incomplete pairwise comparison matrix represented by a complete tripartite graph, where two partitioned sets have one node each, corresponding to the best and worst alternatives, respectively.

A parsimonious approach proposed by \cite{AbastanteCorrenteGrecoIsizakaLami2019} combines direct rating and pairwise comparisons to certain alternatives identified as reference. Therefore, the potential distortions coming from comparing more relevant alternatives with less relevant alternatives are avoided.
Heuristic rating estimation \cite{Kulakowski2015b} assumes \emph{a priori} known and fixed priorities for some alternatives, while the weights of the other alternatives need to be estimated based on pairwise comparisons.

\section{Completion and weighting methods} \label{Sec3}

Analogous to (complete) pairwise comparison matrices, incomplete pairwise comparison matrices are mainly used to obtain a weight vector that reflects the importance of the alternatives on a common scale.
This can be achieved in two different ways:
(a) directly deriving the priorities from the incomplete pairwise comparisons;
(b) completing the matrix first, followed by applying a standard weighting method.

\cite{Harker1987a} has generalised the eigenvector method to the incomplete case as follows. The weight vector is the principal eigenvector of matrix $\mathbf{H} = \left[ h_{ij} \right]$, which has the same size as the incomplete pairwise comparison matrix $\mathbf{A}$.
$h_{ii}$ equals one plus the number of missing elements in the $i$th row of $\mathbf{A}$, $h_{ij} = 0$ for all missing elements $a_{ij} = \ast$, and $h_{ij} = a_{ij}$ for all non-diagonal ($i \neq j$) and known pairwise comparisons in $\mathbf{A}$. For instance, Example \ref{Examp1} implies the following matrix:
\[
\mathbf{H} = \left[
\begin{array}{K{3em} K{3em} K{3em} K{3em}}
    2    	& a_{12}  	& 0  	& a_{14} \\
    a_{21}	& 2       	& a_{23}	& 0   \\
   0		& a_{32}	& 2      	& a_{34} \\
    a_{41}	& 0  	& a_{43}	& 2 \\
\end{array}
\right].
\]

In the complete case, one of the most popular weighting procedures is the logarithmic least squares method \cite{CrawfordWilliams1985, DeGraan1980, deJong1984, Rabinowitz1976, WilliamsCrawford1980}. It minimises the squared differences of the logarithms of matrix entries and their approximations, and has a natural extension by restricting the objective function to the set of known comparisons \cite{Kwiesielewicz1996, TakedaYu1995}. Consequently, the weight vector $\mathbf{w} = \left[ w_i \right]$ is the optimal solution to the following problem:
\begin{align} \label{eq_LLSM}
\min & \sum_{i,j: \, a_{ij} \neq \ast} \left[ \log a_{ij} - \log \left( \frac{w_i}{w_j} \right) \right]^2 \nonumber \\
\text{subject to } & w_i > 0 \text{ for all } i=1,2, \dots n.
\end{align}

The incomplete logarithmic least squares method is the natural extension of the logarithmic least squares method. \cite{LundySirajGreco2017} has shown that the geometric mean of weight vectors, calculated from all spanning trees of the complete graph representing a (complete) pairwise comparison matrix, a method proposed by \cite{Tsyganok2000,Tsyganok2010}, is the unique optimal solution to the logarithmic least squares optimisation problem. According to \cite{BozokiTsyganok2019}, this equivalence holds even for incomplete pairwise comparison matrices: the optimal solution to \eqref{eq_LLSM} is the geometric mean of the weight vectors determined by the spanning trees of graph $G$ representing the incomplete pairwise comparison matrix. 

The logarithmic least squares optimal completion is the pairwise comparison matrix $\mathbf{B} = \left[ b_{ij} \right]$ such that $b_{ij} = a_{ij}$ if $a_{ij} \neq \ast$ and $b_{ij} = w_i / w_j$ if $a_{ij} = \ast$ with $\mathbf{w} = \left[ w_i \right]$ being the optimal solution to \eqref{eq_LLSM}.

Since there exists a one-to-one mapping between weight vectors of order $n$ and consistent pairwise comparison matrices of size $n$, another plausible approach to complete an incomplete pairwise comparison matrix is replacing the $m$ missing entries with positive variables collected in vector $\mathbf{x}$, and minimising the inconsistency of the resulting pairwise comparison matrix $\mathbf{A}(\mathbf{x})$.
Indeed, the logarithmic least squares optimal completion minimises the geometric inconsistency index \cite{CrawfordWilliams1985, AguaronMoreno-Jimenez2003}.

The most popular inconsistency index, proposed by \cite{Saaty1977}, is a monotonic function of the dominant eigenvalue (see Section~\ref{Sec2}). This leads to the following optimisation problem \cite{ShiraishiObata2002, ShiraishiObataDaigo1998}:
\begin{align} \label{eq_EM}
\min_{\mathbf{x} > 0} \lambda_{\max} \left( \mathbf{A}(\mathbf{x}) \right) \nonumber \\
\lambda_{\max} \left( \mathbf{A}(\mathbf{x}) \right) \mathbf{w} = \mathbf{A}(\mathbf{x}) \mathbf{w}.
\end{align}
It is called the incomplete eigenvalue method; the variables in $\mathbf{x}$ minimise the dominant eigenvalue of the associated pairwise comparison matrix, while the weight vector coincides with the (normalised) right eigenvector of this matrix.

Consequently, the eigenvalue optimal completion is the pairwise comparison matrix $\mathbf{B} = \left[ b_{ij} \right]$ such that $b_{ij} = a_{ij}$ if $a_{ij} \neq \ast$ and $b_{ij}$ equals the appropriate entry of vector $\mathbf{x}$, the optimal solution to \eqref{eq_EM} if $a_{ij} = \ast$.

Both problems~\eqref{eq_LLSM} and \eqref{eq_EM} have a unique solution if and only if the graph associated with the incomplete pairwise comparison matrix is connected \cite[Theorems~2 and 4]{BozokiFulopRonyai2010}. This is a reasonable condition: any two alternatives can be compared at least indirectly, through other alternatives. The necessary and sufficient condition for the uniqueness of the solution to problem~\eqref{eq_LLSM} has already been proved by \cite{KaiserSerlin1978}.

The logarithmic least squares optimal completion can be obtained by solving a system of linear equations, see \cite{KaiserSerlin1978} and \cite[Theorem~4]{BozokiFulopRonyai2010}. \cite{BozokiFulopRonyai2010} provides an algorithm for finding the eigenvalue optimal completion based on the method of cyclic coordinates.

Recent studies \cite{TuWuPedrycz2023, XuWang2024} argue that the logarithmic least squares method is a promising technique to derive weights from a best--worst method matrix due to
(a) its simple calculation as the solution of a system of linear equations;
(b) the uniqueness of the weights;
(c) accounting for indirect comparisons.

The relation between the logarithmic least squares and the eigenvalue optimal completions is also worth investigating.
First, take the case $n=3$. If one comparison is unknown, but the representing graph remains connected, then it should be a spanning tree, and there exists a unique consistent completion $\mathbf{B} = \left[ b_{ij} \right]$: $a_{ik} = \ast$ implies $b_{ik} = a_{ij} a_{jk}$. Thus, the optima of \eqref{eq_LLSM} and \eqref{eq_EM}, zero and $n$, respectively; each attains its theoretical minimums.
On the other hand, if two comparisons are missing, the representing graph is not connected, and infinitely many consistent completion exists.

Surprisingly, the two optimal completions coincide even if $n=4$, independently of the number of unknown comparisons \cite[Theorem~1]{CsatoAgostonBozoki2024}.
However, they might be different for at least five alternatives with one missing comparison only.

\begin{example} \label{Examp2}
\cite[Lemma~2]{CsatoAgostonBozoki2024}
Take the following incomplete pairwise comparison matrix:
\[
\mathbf{A} = \left[
\begin{array}{K{3em} K{3em} K{3em} K{3em} K{3em}}
    1     &  1/2  & 5     &  1/6  & \ast \\
    2     & 1     & 4     &  1/2  &  1/6 \\
     1/5  &  1/4  & 1     &  1/6  &  1/7 \\
    6     & 2     & 6     & 1     &  1/2 \\
    \ast  & 6     & 7     & 2     & 1     \\
\end{array}
\right].
\]
Denote by $\mathbf{B}$ and $\mathbf{C}$ the logarithmic least squares and the eigenvalue optimal completions, respectively. Solving problems~\eqref{eq_LLSM} and \eqref{eq_EM} leads to $b_{15} = 0.1705$ and $c_{15} = 0.1798$.
\end{example}

The definition of inconsistency involves three alternatives. Therefore, several inconsistency indices are based on triads, pairwise comparison (sub)matrices determined by three alternatives. For example, the Koczkodaj inconsistency index \cite{DuszakKoczkodaj1994, Koczkodaj1993}---which has been axiomatically characterised by \cite{Csato2018a}---depends on the most inconsistent triad.
\cite{Csato2019c} shows that there exists only one natural inconsistency ranking on the set of triads, given by $\mathit{TI} = \max \left\{ a_{ik} / \left( a_{ij} a_{jk} \right); \left( a_{ij} a_{jk} \right) / a_{ik} \right\}$, and presents its axiomatic characterisation. Unsurprisingly, almost all inconsistency indices are functionally related to $\mathit{TI}$ on the set of triads \cite{Cavallo2020}.

The uniqueness of triad inconsistency ranking suggests another plausible procedure to determine the missing comparisons.
The \emph{lexicographically optimal completion} \cite{AgostonCsato2024}, inspired by the solution concept nucleolus in cooperative game theory \cite{Schmeidler1969}, minimises the inconsistency of the most inconsistent triad first, followed by the inconsistency of the second most inconsistent triad, and so on.

Formally, let $\mathit{TI}_{ijk}(\mathbf{x})$ be the inconsistency of the triad determined by the three alternatives $1 \leq i,j,k \leq n$ according to the measure $\mathit{TI}$, as a function of the missing comparisons collected in vector $\mathbf{x}$.
Denote the vector of the $n(n-1)(n-2)/6$ local inconsistencies, arranged in non-increasing order, by $\theta(\mathbf{x})$. Thus, $u < v$ implies $\theta_u(\mathbf{x}) \geq \theta_v(\mathbf{x})$.
Matrix $\mathbf{B} = \mathbf{A}(\mathbf{x})$ is called a \emph{lexicographically optimal completion} of the incomplete pairwise comparison matrix $\mathbf{A}$ if, for any other completion $\mathbf{A}(\mathbf{y})$, $\theta_u(\mathbf{x}) = \theta_u(\mathbf{y})$ for all $u < v$ implies $\theta_v(\mathbf{x}) \leq \theta_v(\mathbf{y})$.

\begin{example} \label{Examp3}
\cite[Examples~2 and 3]{AgostonCsato2024}
Consider the following incomplete pairwise comparison matrix, where the two missing entries above the diagonal are substituted by variables:
\[
\mathbf{A(\mathbf{x})} = \left[
\begin{array}{K{3em} K{3em} K{3em} K{3em}}
    1     	& 2			& x_{13}   	& x_{14} \\
    1/2		& 1       	& 1			& 8 \\
   1/x_{13}	& 1 		& 1      	& 1 \\
   1/x_{14}	& 1/8	 	& 1  		& 1 \\
\end{array}
\right].
\]
This matrix contains four triads with the following inconsistencies:
\[
\mathit{TI}_{123}(\mathbf{x}) = \max \left\{ \frac{x_{13}}{2}; \frac{2}{x_{13}} \right\};
\]
\[
\mathit{TI}_{124}(\mathbf{x}) = \max \left\{ \frac{x_{14}}{16}; \frac{16}{x_{14}} \right\};
\]
\[
\mathit{TI}_{134}(\mathbf{x}) = \max \left\{ \frac{x_{14}}{x_{13}}; \frac{x_{13}}{x_{14}} \right\};
\]
\[
\mathit{TI}_{234}(\mathbf{x}) = \max \left\{ 8; \frac{1}{8} \right\} = 8.
\]

The unique lexicographically optimal filling is $x_{13} = 4$ and $x_{14} = 8$ with the vector of inconsistencies $\theta(\mathbf{x}) = \left[ 8, 2, 2, 2 \right]$.
First, note that $\theta_1(\mathbf{y})$ cannot be smaller than 8 due to $\mathit{TI}_{234}(\mathbf{x})$, independently of the unknown comparisons $\mathbf{x}$.
Second, $x_{13} > 4$ would lexicographically increase $\theta(\mathbf{x})$ since $\mathit{TI}_{123}(\mathbf{x}) > 2$.
Analogously, $x_{14} < 8$ would lexicographically increase $\theta(\mathbf{x})$ since $\mathit{TI}_{124}(\mathbf{x}) > 2$.
Finally, $x_{13} \leq 4$ and $x_{14} \geq 8$ imply $\mathit{TI}_{134}(\mathbf{x}) \geq 2$, and the inequality is strict if $x_{13} < 4$ or $x_{14} > 8$.
\end{example}

The lexicographically optimal completion can be computed by solving successive linear programming (LP) problems if the logarithmically transformed entries of the original pairwise comparison matrix are used \cite[Section~3]{AgostonCsato2024}:
\begin{enumerate}
\item \label{LP_solution1}
A linear programming problem is solved to minimise the natural triad inconsistency index for all triads with a still unknown value of $\mathit{TI}$.
\item
A triad, whose inconsistency index $\mathit{TI}$ cannot be smaller than the optimum of the current LP, which is indicated by the non-zero shadow price of at least one corresponding constraint, is chosen.
\item
The inconsistency index $\mathit{TI}$ is fixed for this triad (or one of these triads), the two associated constraints are removed from the LP, and we return to Step~\ref{LP_solution1}.
\end{enumerate}
Since the number of constraints continuously decreases, this procedure requires solving at most as many LPs as the number of triads, $n(n-1)(n-2)/6$. The running time of the algorithm remains below one second for most problems arising in practice.
In addition, the lexicographically optimal completion can be calculated without using an LP solver if the missing entries are independent (occur in different rows and columns), as in Example~\ref{Examp3} \cite[Section~5]{AgostonCsato2024}.

Naturally, several further completion methods have been proposed in the literature. \cite{TekileBrunelliFedrizzi2023} compares eleven completion methods, including the incomplete eigenvector and the incomplete logarithmic least squares methods, based on numerical simulations and hierarchical clustering.

\section{The issue of ordinal violations} \label{Sec4}

Let $\mathbf{A} = \left[ a_{ij} \right]$ be an incomplete pairwise comparison matrix.
Its weight vector $\mathbf{w} = \left[ w_i \right]$ shows an \emph{ordinal violation} if there exist alternatives $i,j$ such that $a_{ij} > 1$ but $w_i \leq w_j$.

Ordinal violations are important because the decision-maker may be unwilling to accept a weight vector if it contradicts the direction of the preferences. While a weight vector without an ordinal violation does not necessarily exist, ordinal violation is particularly problematic in at least two cases:
(a) when the pairwise comparisons are given by a directed acyclic graph (Section~\ref{Sec41}) \cite{CsatoRonyai2016};
(b) in a best--worst method matrix (Section~\ref{Sec42}) \cite{TuWuPedrycz2023}.

\subsection{Incomplete pairwise comparisons implied by connected directed acyclic graphs} \label{Sec41}

A \emph{directed graph} is a tuple $(N,E)$, where $N$ is the set of vertices, and $E$ is the set of ordered pairs of vertices (arcs).
A \emph{cycle} is a sequence of arcs in $E$ such that
(i) the ending vertex of any arc in the sequence coincides with the starting vertex of the next arc; and
(ii) the starting vertex of the first arc is the ending vertex of the last arc.
A directed graph $(N,E)$ is \emph{acyclic} if it does not contain any cycle.
A directed graph $(N,E)$ is (weakly) \emph{connected} if the corresponding undirected graph, where the directions of all arcs are deleted, is connected.

Any connected directed acyclic graph (CDAG) generates an incomplete pairwise comparison matrix, with parameter $\alpha > 1$ reflecting the dominance relation determined by the arcs.
Formally, let $(N,E)$ be a CDAG with $|N| = n$ vertices.
The associated \emph{CDAG-based incomplete pairwise comparison matrix} is $\mathbf{A}(\alpha) = \left[ a_{ij} \right]$ such that, for all $1 \leq i,j \leq n$:
\begin{itemize}
\item
$a_{ii} = 1$;
\item
$a_{ij} = \alpha$ and $a_{ji} = 1 / \alpha$ if $(i,j) \in E$;
\item
$a_{ij} = \ast$ and $a_{ji} = \ast$ if $(i,j) \notin E$ and $(j,i) \notin E$.
\end{itemize}

In a connected directed acyclic graph, it can be assumed without loss of generality that $(i,j) \in E$ implies $i < j$.
Then the associated CDAG-based incomplete pairwise comparison matrix has entries higher (lower) than one only in the upper (lower) triangle.

The weight vector $\mathbf{w}^{\mathit{LLSM}}$ given by the incomplete logarithmic least squares method~\eqref{eq_LLSM} may show an ordinal violation when it is applied to a CDAG-based incomplete pairwise comparison matrix.

\begin{figure}[t!]
\centering

\begin{tikzpicture}[scale=1, auto=center, transform shape, >=triangle 45]
\tikzstyle{every node}=[draw,shape=circle];
  \node (n1)  at (2*360/7:3) {$1$};
  \node (n2)  at (360/7:3)   {$2$};
  \node (n3)  at (0:3)   	 {$3$};
  \node (n4)  at (6*360/7:3) {$4$};
  \node (n5)  at (5*360/7:3) {$5$};
  \node (n6)  at (4*360/7:3) {$6$};
  \node (n7)  at (3*360/7:3) {$7$};

  \foreach \from/\to in {n1/n2,n1/n6,n1/n7,n2/n3,n2/n4,n3/n4,n3/n5,n4/n5,n4/n6,n5/n6,n5/n7}
    \draw [->] (\from) -- (\to);
\end{tikzpicture}

\caption{The directed acyclic graph of Example~\ref{Examp4}}
\label{Fig2}
\end{figure}
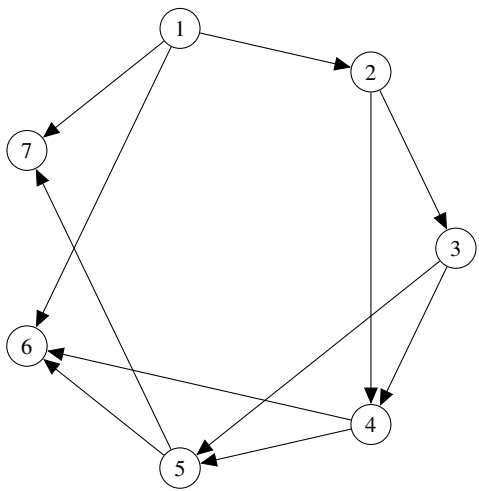

\begin{example} \label{Examp4}
\cite[Example~3.4]{CsatoRonyai2016}
Consider the directed acyclic graph in Figure~\ref{Fig2}.
The associated CDAG-based incomplete pairwise comparison matrix is:
\[
\mathbf{A}(\alpha) =
\left(
\begin{array}{K{3em}K{3em}K{3em}K{3em}K{3em}K{3em}K{3em}}
    1     & \alpha     & \ast & \ast & \ast & \alpha     & \alpha \\
    1/\alpha   & 1     & \alpha     & \alpha     & \ast & \ast & \ast \\
    \ast & 1/\alpha   & 1     & \alpha     & \alpha     & \ast & \ast \\
    \ast & 1/\alpha   & 1/\alpha   & 1     & \alpha     & \alpha     & \ast \\
    \ast & \ast & 1/\alpha   & 1/\alpha   & 1     & \alpha     & \alpha \\
    1/\alpha   & \ast & \ast & 1/\alpha   & 1/\alpha   & 1     & \ast \\
    1/\alpha   & \ast & \ast & \ast & 1/\alpha   & \ast & 1 \\
\end{array}
\right).
\]
The logarithms of the weights given by the incomplete logarithmic least squares method are
\[
\mathbf{y}(\mathbf{A}) = \left[
\begin{array}{K{2em}K{2em}K{2em}K{2em}K{2em}K{2em}K{2em}}
    34    & 36    & 24    & 1     & -14   & -42   & -39 \\
\end{array}
\right]^\top \log \alpha / 49.
\]
Hence, $w_1^{\mathit{LLSM}}(\mathbf{A}) < w_2^{\mathit{LLSM}}(\mathbf{A})$, which is an ordinal violation due to $a_{12} = \alpha > 1$.
\end{example}

Example~\ref{Examp4} is minimal with respect to the number of alternatives and, among them, with respect to the number of known comparisons (11). There exist analogous examples with eight alternatives but only 10 known comparisons, see \cite[Example~3.7]{CsatoRonyai2016}.

If the incomplete logarithmic least squares method is used, the ranking of the alternatives does \emph{not} depend on the value of parameter $\alpha > 1$ \cite[Corollary~3.2]{CsatoRonyai2016}.
Therefore, the ranking of the alternatives by the logarithmic least squares method is unique for any CDAG-based incomplete pairwise comparison matrix.

The weight vector $\mathbf{w}^{\mathit{EM}}$ given by the incomplete eigenvector method~\eqref{eq_EM} may also exhibit an ordinal violation.

\begin{figure}[t!]
\centering
\caption{The directed acyclic graph of Example~\ref{Examp5}}
\label{Fig3}
\begin{tikzpicture}[scale=1, auto=center, transform shape, >=triangle 45]
\tikzstyle{every node}=[draw,shape=circle];
  \node (n1)  at (112.5:3)  {$1$};
  \node (n2)  at (67.5:3)   {$2$};
  \node (n3)  at (22.5:3)   {$3$};
  \node (n4)  at (337.5:3)  {$4$};
  \node (n5)  at (292.5:3)  {$5$};
  \node (n6)  at (247.5:3)  {$6$};
  \node (n7)  at (202.5:3)  {$7$};
  \node (n8)  at (157.5:3)  {$8$};

  \foreach \from/\to in {n1/n2,n1/n7,n1/n8,n2/n3,n2/n4,n3/n5,n3/n6,n4/n5,n4/n6,n5/n7,n5/n8,n6/n7,n6/n8}
    \draw [->] (\from) -- (\to);
\end{tikzpicture}
\end{figure}
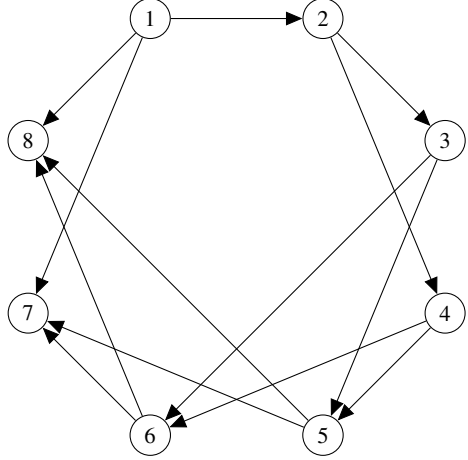

\begin{example} \label{Examp5}
\cite[Theorem~4.2 and Proposition 4.3]{CsatoRonyai2016}
Consider the directed acyclic graph in Figure~\ref{Fig3}.

Solving optimisation problem~\eqref{eq_EM} results in
\[
\mathbf{w}^{EM} \left( \mathbf{A} (\alpha = 3) \right) = \left[
\begin{array}{K{2.5em}K{2.5em} K{2.5em}K{2.5em} K{2.5em}K{2.5em} K{2.5em}K{2.5em}}
    24.04 & 24.42 & 14.81 & 14.81 & 7.29 & 7.29 & 3.67 & 3.67 \\
\end{array}
\right]^\top,
\]
and
\[
\mathbf{w}^{EM} \left( \mathbf{A} (\alpha = 4) \right) = \left[
\begin{array}{K{2.5em}K{2.5em} K{2.5em}K{2.5em} K{2.5em}K{2.5em} K{2.5em}K{2.5em}}
    28.28 & 26.56 & 14.04 & 14.04 & 5.94 & 5.94 & 2.60 & 2.60 \\
\end{array}
\right]^\top,
\]
where the sum of the weights is normalised to 100.
Consequently, $w_1^{EM} \left( \mathbf{A} (\alpha = 3) \right) < w_2^{EM} \left( \mathbf{A} (\alpha = 3) \right)$ but $a_{12} = 3$, which implies an ordinal violation.
\end{example}

In Example~\ref{Examp5}, $w_1^{EM} \left( \mathbf{A} (\alpha = 3) \right) < w_2^{EM} \left( \mathbf{A} (\alpha = 3) \right)$ and $w_1^{EM} \left( \mathbf{A} (\alpha = 4) \right) > w_2^{EM} \left( \mathbf{A} (\alpha = 4) \right)$. Hence, the ranking of the alternatives by the eigenvector method is not necessarily unique for a CDAG-based incomplete pairwise comparison matrix.

Examples~\ref{Examp4} and \ref{Examp5} suggest that neither the logarithmic least squares method, nor the eigenvector method is ideal for deriving weights from a CDAG-based incomplete pairwise comparison matrix.
However, the weight vector is guaranteed to be free of any ordinal violation if the lexicographically optimal completion is combined with either the eigenvector or the logarithmic least squares method to derive the priorities from the complete pairwise comparison matrix \cite[Theorem~1]{Csato2024b}. This seems to be a strong argument in favour of using the lexicographically optimal completion.

\subsection{Incomplete pairwise comparisons implied by the best--worst method} \label{Sec42}

Let $\mathbf{A} = \left[ a_{ij} \right]$ be a best--worst method matrix, where the first alternative is the best and the last, $n$th alternative is the worst: $a_{1j} > 1$ for all $2 \leq j \leq n$, and $a_{jn} > 1$ for all $1 \leq j \leq n-1$. Furthermore, $a_{ij} = \ast$ holds for all $2 \leq i,j \leq n-1$ if $i \neq j$.
If the weight vector $\mathbf{w} = \left[ w_i \right]$ shows an \emph{ordinal violation}, that is, $w_1 < w_j$ or $w_n > w_j$ for a particular index $2 \leq j \leq n-1$, then the best (worst) alternative---identified in the first stage---does not receive the highest (lowest) priority based on the best--worst method matrix in the final phase. Such an anomaly certainly calls for revisiting the pairwise comparisons, which can be frustrating and time-consuming in any application.

In the following, the known entries of the best--worst method matrix are assumed to be between 1/9 and 9 as proposed by \cite{Saaty1980}, but they are not necessarily integers or their reciprocals.
\cite{Csato2025a} derives two sufficient conditions to ensure that the logarithmic least squares weights do not show an ordinal violation for a best--worst method matrix.

The first result contains a uniform lower bound for the preference of the best alternative over all other alternatives, as well as for the preference of all other alternatives over the worst alternative. The maximal numerical preference is also constrained by a polynomial of this uniform lower bound.
In particular, if $a_{1j} \geq p$ (the best alternative is at least $p$ times better than all other alternatives) and $a_{jn} \geq p$ (the worst alternative is at least $p$ times worse than all other alternatives), while the maximal preference between two alternatives is at most $p^3$, then the logarithmic least squares weights do not contain any ordinal violation for a best--worst method matrix \cite[Theorem~1]{Csato2025a}.

Besides the uniform lower bound for the intensity of all preferences and the upper bound for the maximal preference, the second result adopts the natural assumption that the preference between the best and the worst alternatives is the strongest.
In particular, if $a_{1j} \geq p$ and $a_{jn} \geq p$, as well as $a_{1n} \geq \max \left\{ a_{1j}: 2 \leq j \leq n-1 \right\}$ and $a_{1n} \geq \max \left\{ a_{jn}: 2 \leq j \leq n-1 \right\}$, but $a_{1n} \leq p^{4/(n-3)+3}$, then the logarithmic least squares weights do not contain any ordinal violation for a best--worst method matrix \cite[Theorem~2]{Csato2025a}.

The proof of both statements is based on a careful consideration of the logarithmic least squares weights, for which a closed-form formula exists \cite[Theorem~4]{BozokiFulopRonyai2010}. Since $p^{4/(n-3)+3} > p^3$, the upper bound is less restrictive in the second theorem, although it contains the additional constraint that the preference between the best and the worst alternatives should be the strongest.
However, $4/(n-3) + 3 \to 3$ when $n \to \infty$; hence, the sufficient condition of \cite[Theorem~2]{Csato2025a} is more restrictive than the sufficient condition of \cite[Theorem~1]{Csato2025a} for a high number of alternatives, due to the additional constraint on the preference between the best and the worst alternatives.
On the other hand, the second result is quite powerful if $n$ is (relatively) small. 
For instance, if the entries of a best--worst method matrix belong to the Saaty scale (integers from 2 to 9 and their reciprocals), then $p \geq 2$ holds for the uniform lower bound, and $2^{4/(n-3) + 3} > 9$ is satisfied unless the number of alternatives $n$ exceeds 26.

\begin{example} \label{Examp6}
\cite[Example~1]{Csato2025a}
Consider the following best--worst method matrix with six alternatives:
\[
\mathbf{A} = \left[
\begin{array}{K{3em} K{3em} K{3em} K{3em} K{3em} K{3em}}
    1     	 & 2	 	& 2    & 2 	  & 2    & 2 \\
    1/2		 & 1     	& \ast & \ast & \ast & 9 \\
    1/2      & \ast    	& 1    & \ast & \ast & 2 \\
    1/2      & \ast    	& \ast & 1    & \ast & 2 \\
    1/2      & \ast  	& \ast & \ast & 1    & 2 \\
    1/2      & 1/9    	& 1/2  & 1/2  & 1/2  & 1 \\
\end{array}
\right].
\]
The logarithmic least squares priorities are
\[
\mathbf{w} = \left[
\begin{array}{K{3em} K{3em} K{3em} K{3em} K{3em} K{3em}}
    26.45 & 27.78	& 13.10 & 13.10 & 13.10 & 6.48 \\
\end{array}
\right]^\top
\]
if their sum is equal to 100. Hence, the second alternative has the highest weight, even though the first has been identified as the best.
\end{example}

Naturally, Example~\ref{Examp6} shows an extreme situation.
Assume that the number of alternatives is six and the pairwise comparisons belong to the Saaty scale of
\[
\{ 1/9, 1/8, \dots ,1/2,1,2, \dots ,8,9 \}.
\]
The number of these best--worst method matrices is $8^9 = 134{,}217{,}728$. Among them, $7^9 = 40{,}353{,}607$ matrices (30.1\%) satisfy the sufficient conditions of \cite[Theorem~1]{Csato2025a} as all nine comparisons can be at most $2^3 = 8$. However, the logarithmic least squares weights show an ordinal violation for only 56 matrices.

\section{Generalised inconsistency thresholds} \label{Sec5}

\begin{table}[t!]
\caption{The random index $\mathit{RI}_n$ for complete pairwise comparison matrices}
\centering
\label{Table1}
    \begin{tabularx}{\textwidth}{l CCCCCCC} \toprule
    Matrix size $n$ & 4     & 5     & 6     & 7     & 8     & 9     & 10 \\ \midrule
    Random index $RI_n$ & 0.884 & 1.109 & 1.249 & 1.341 & 1.404 & 1.451 & 1.486 \\ \bottomrule
    \end{tabularx}
\end{table}

As we have seen in Section~\ref{Sec2}, the acceptance threshold 0.1 applies to the inconsistency ratio $\mathit{CR}$, which depends on the random index $\mathit{RI}$, the average inconsistency index $\mathit{CI}$ of a large number of randomly generated pairwise comparison matrices. 
However, since the incomplete eigenvalue method chooses the missing entries such that the dominant eigenvalue and, consequently, the value of $\mathit{CI}$ is minimised (see Section~\ref{Sec3}, in particular, optimisation problem~\eqref{eq_EM}), the random index $\mathit{RI}_n$, computed from complete pairwise comparison matrices, cannot be used for an incomplete pairwise comparison matrix. Therefore, by adopting the values reported in Table~\ref{Table1}, the ratio of incomplete pairwise comparison matrices with an acceptable level of inconsistency will exceed the original intention of Saaty---and this deviation increases with the number of unknown comparisons. For example, when the representing graph $G$ is a spanning tree with $n$ nodes and $n-1$ edges, the corresponding incomplete matrix has a consistent completion, and the random index $\mathit{RI}$ should be equal to zero.

Hence, the random index needs to be recomputed for incomplete pairwise comparison matrices, which is done by \cite{AgostonCsato2022}, a paper that has received the Best Paper Award in the journal \emph{Omega} in 2022.
There are three different ways to choose the missing entries. \cite{AgostonCsato2022} restricts their value to the interval $\left[ 1/9, 9 \right]$, but they should not be integers or their reciprocals. Thus, the value of the random index remains positive even if the incomplete pairwise comparison matrix is represented by a spanning tree: although these matrices always have a unique consistent completion, it may be infeasible due to the lower and upper bounds on the missing entries.

\begin{table}[t!]
\caption{The random index for incomplete pairwise comparison matrices \cite[Table~2]{AgostonCsato2022}}
\centering
\label{Table2}
\begin{threeparttable}
\rowcolors{3}{gray!20}{}
    \begin{tabularx}{\textwidth}{c CCCC} \toprule \hiderowcolors    
    \multirow{2}{*}{Missing elements $m$} & \multicolumn{4}{c}{Matrix size $n$} \\
          & 4     & 5     & 6     & 7     \\ \bottomrule \showrowcolors
    0     & 0.884 & 1.109 & 1.249 & 1.341 \\
    1     & 0.583 (0.531) & 0.925 (0.485) & 1.128 (0.400) & 1.256 (0.330) \\
    2     & 0.306 (0.387) & 0.739 (0.452) & 1.007 (0.392) & --- \\
    3     & 0.053 (0.073) & 0.557 (0.405) & 0.883 (0.380) & --- \\
    4     & ---   & 0.379 (0.340) & 0.758 (0.364) & --- \\
    5     & ---   & 0.212 (0.247) & 0.634 (0.344) & --- \\
    6     & ---   & 0.059 (0.068) & 0.510 (0.317) & --- \\
    7     & ---   & ---   & 0.389 (0.281) & --- \\
    8     & ---   & ---   & 0.271 (0.234) & --- \\
    9     & ---   & ---   & 0.161 (0.170) & --- \\ \toprule
    \end{tabularx}
\begin{tablenotes} \footnotesize
\item
All values are based on 1 million randomly generated matrices. Standard deviations are in parenthesis.
\end{tablenotes}
\end{threeparttable}
\end{table}

The generalised random indices are given in Table~\ref{Table2} as a function of the number of alternatives $n$ and the number of missing entries $m$.
Both factors have a non-negligible effect; the values of $RI_{n,m}$ are statistically different from each other at all commonly used levels of significance.

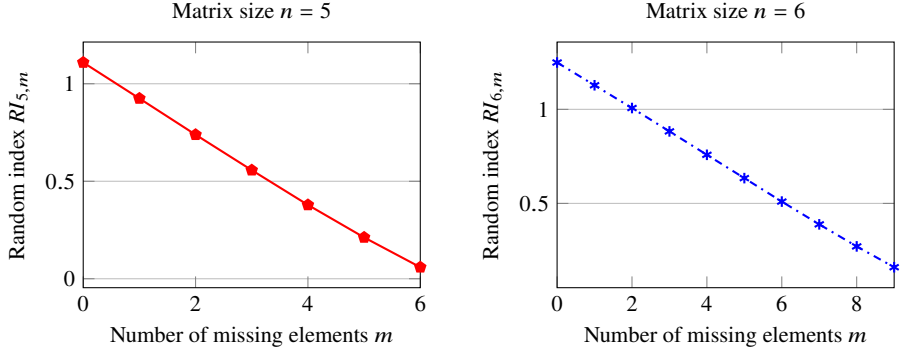
\begin{figure}[t!]

\begin{tikzpicture}
\begin{axis}[
name = axis1,
title = {Matrix size $n=5$},
title style = {font=\small},
xlabel = Number of missing elements $m$,
x label style = {font=\small},
ylabel = Random index $\mathit{RI}_{5,m}$,
y label style = {font=\small},
width = 0.5\textwidth,
height = 0.4\textwidth,
ymajorgrids = true,
xmin = 0,
xmax = 6,
%ymin = 0,
%ymax = 100,
%max space between ticks=50,
%legend style = {font=\small,at={(0.05,-0.15)},anchor=north west,legend columns=8},
%legend entries = {$S1 \quad$,$S2 \quad$,$S3 \quad$,$S4 \quad$,$F \quad$,$G1 \quad$,$G2 \quad$,$G3$}
] 
% 5 x 5
\addplot [red, thick, mark=pentagon*] coordinates {
(0,1.109)
(1,0.9246087)
(2,0.7387896)
(3,0.556883)
(4,0.37861)
(5,0.211943)
(6,0.0591422)
};
\end{axis}

\begin{axis}[
at = {(axis1.south east)},
xshift = 0.15\textwidth,
title = {Matrix size $n=6$},
title style = {font=\small},
xlabel = Number of missing elements $m$,
x label style = {font=\small},
ylabel = Random index $\mathit{RI}_{6,m}$,
y label style = {font=\small},
width = 0.5\textwidth,
height = 0.4\textwidth,
ymajorgrids = true,
xmin = 0,
xmax = 9,
]
% 6 x 6
\addplot [blue, thick, dashdotted, mark=asterisk, mark options={solid,thick}] coordinates {
(0,1.249)
(1,1.127994)
(2,1.007017)
(3,0.8827498)
(4,0.7582752)
(5,0.6343154)
(6,0.50963)
(7,0.388586)
(8,0.2713317)
(9,0.160653)
};
\end{axis}
\end{tikzpicture}

%\captionsetup{justification=centering}
\caption{The relationship between the random index and the number of missing entries}
\label{Fig4}
\end{figure}

According to Figure~\ref{Fig4}, the random index is an almost linear, monotonically decreasing function of $m$ for both five and six alternatives. This observation gives a relatively good approximation for the random index $\mathit{RI}_{n,m}$ based on the well-known values in Table~\ref{Table1} \cite[Table~3]{AgostonCsato2022}:
\begin{equation} \label{RI_approx}
RI_{n,m} \approx \left[ 1 - \frac{2m}{(n-1)(n-2)} \right] RI_{n,0}.
\end{equation}
Note that formula~\eqref{RI_approx} gives zero if the number of unknown comparisons is $m = (n-1)(n-2)/2$ (when the representing graph is a spanning tree since it is assumed to be connected), while \eqref{RI_approx} equals $\mathit{RI}_n$ in the absence of missing entries ($m=0$).

Approximation~\eqref{RI_approx} ensures that the generalised inconsistency thresholds can be used in almost all applications as $\mathit{RI}_n$ has been computed for $n \leq 16$ \cite{AguaronMoreno-Jimenez2003} and for $n \leq 15$ \cite{AlonsoLamata2006}. In practice, the optimal number of alternatives does rarely exceed nine \cite{SaatyOzdemir2003}.

One of the most important roles of these generalised inconsistency thresholds is allowing continuous monitoring of inconsistency.
\cite{BozokiDezsoPoeszTemesi2013} conducted an experiment, where we know not only the final complete pairwise comparison matrices, but all incomplete submatrices after a given number of answers. 

\begin{example} \label{Examp7}
Consider the following pairwise comparison matrix, which collects the numerical answers on how much more a summer house is favoured against another summer house:
\[
\mathbf{A} = \left[
\begin{array}{K{3em} K{3em} K{3em} K{3em} K{3em} K{3em}}
    1     & 2     & 7     & 7     & 7     & 5 \\
     1/2  & 1     & 5     & 7     & 5     & 2 \\
     1/7  &  1/5  & 1     & 1     &  1/5  &  1/5 \\
     1/7  &  1/7  & 1     & 1     &  1/3  &  1/3 \\
     1/7  &  1/5  & 5     & 3     & 1     & 3 \\
     1/5  &  1/2  & 5     & 3     &  1/3  & 1 \\
\end{array}
\right].
\]
The pairwise comparisons have been obtained by the method of \cite{Ross1934} in order to ensure a balanced appearance of any alternative throughout the collection process, as well as to ensure that they occupy both the first and the second positions in the pairwise comparisons roughly with an equal probability. Thus, the implied questioning order is $a_{12}$, $a_{64}$, $a_{51}$, $a_{32}$, $a_{56}$, $a_{13}$, $a_{24}$, $a_{61}$, $a_{43}$, $a_{52}$, $a_{14}$, $a_{35}$, $a_{26}$, $a_{45}$, and $a_{36}$.

\begin{figure}[t!]
\centering
\begin{tikzpicture}
\begin{axis}[
name = axis1,
title style = {font=\small},
xlabel = Number of known entries $n(n-1)/2 - m$,
x label style = {font=\small},
width = \textwidth,
height = 0.6\textwidth,
ymajorgrids = true,
xmin = 6,
xmax = 15,
ymin = 0,
%ymax = 100,
%max space between ticks=50,
yticklabel style = {scaled ticks=false,/pgf/number format/fixed,/pgf/number format/precision=2}, 
legend style = {font=\small,at={(0.05,-0.15)},anchor=north west,legend columns=1},
legend entries = {$\mathit{CR}$ without accounting for the number of missing entries based on $\mathit{RI}_n$,$\mathit{CR}$ with accounting for the number of missing entries based on $\mathit{RI}_{n,m} \;$},
] 
% original CI
\addplot [red, thick, mark=pentagon*] coordinates {
(6,0.00569675820656525)
(7,0.0392998983186549)
(8,0.048054261008807)
(9,0.0639635228182546)
(10,0.0661744251401121)
(11,0.0669433706965572)
(12,0.0810705892714171)
(13,0.0894788678943154)
(14,0.0913774395516413)
(15,0.093606)
};
% Our CI
\addplot [blue, thick, dashdotted, mark=asterisk, mark options={solid,thick}] coordinates {
(6,0.0441941)
(7,0.1811276)
(8,0.1542925)
(9,0.1566479)
(10,0.1303657)
(11,0.1103064)
(12,0.114674)
(13,0.1109822)
(14,0.1011795)
(15,0.093606)
};
\end{axis}
\end{tikzpicture}

%\captionsetup{justification=centering}
\caption{The evolution of the inconsistency ratio $\mathit{CR}$ in Example~\ref{Examp7}}
\label{Fig5}
\end{figure}
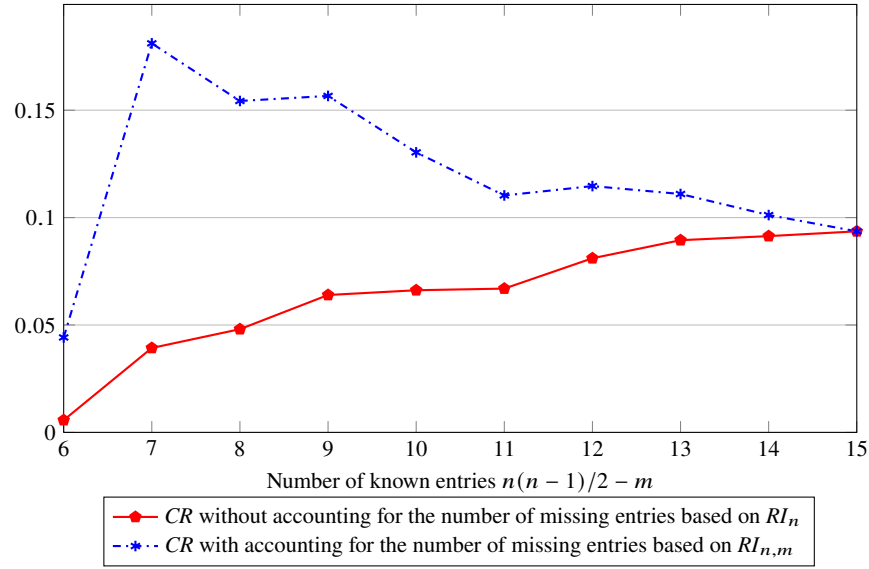

The value of the inconsistency ratio $\mathit{CR}$ during the filling-in procedure is plotted in Figure~\ref{Fig5}.
The solid red line is based on the na\"ive approach \cite[Figure~2]{BozokiDezsoPoeszTemesi2013}, when the random index from Table~\ref{Table1}---which is valid for complete pairwise comparison matrices---is not adjusted according to the number of missing entries $m$. Under the eigenvalue-optimal completion, it is guaranteed to monotonically increase, because the decision-maker cannot choose the value of the next pairwise comparison such that the minimal value of the dominant eigenvalue (consequently, the inconsistency index $\mathit{CI}$ and the inconsistency ratio $\mathit{CR}$) remains below the optimum of minimisation problem~\eqref{eq_EM}. Hence, the inconsistency ratio appears to remain below the threshold of 0.1 throughout the experiment.

On the other hand, the dashed blue line uses the generalised random indices $\mathit{RI}_{n,m}$ from Table~\ref{Table2}. Then, a sudden jump of inconsistency emerges after the decision-maker carries out the seventh comparison ($a_{24}$). Although the inconsistency of the final, complete pairwise comparison matrix can be tolerated, continuous monitoring draws our attention to the fact that the value of $a_{24}$ is worth reconsidering.
\end{example}

A recent working paper investigates how the positions of the missing entries influence the value of the random index for incomplete pairwise comparison matrices \cite{AgostonCsato2026}. While their impact is smaller compared to the number of missing entries, the graph structure still significantly affects the ratio of randomly generated pairwise comparison matrices with an acceptable level of inconsistency.

Further inconsistency indices have also been extended for incomplete pairwise comparisons \cite{KulakowskiTalaga2020}, but they usually do not have such well-established thresholds.

\section{Optimal filling-in sequences} \label{Sec6}

\input{Figure6_optimal_filling_in_sequences}

The filling-in process of a pairwise comparison matrix is of theoretical and practical interest. How should one balance between the time of the decision maker and the quality of information? 
Furthermore, some (online) questionnaires are not always filled in completely because the respondent quits earlier. Preparing for such cases, \cite{SzadoczkiBozoki2025} investigates the filling-in process up to six alternatives with the help of a graph of graphs approach.
Optimal subsets of comparisons---that result, on average, in the closest weight vectors to the complete case---are identified, and optimal sequences, including as many optimal subsets as possible, are suggested. Figure~\ref{Fig6} presents the optimal filling-in patterns for the case of $n=5$.

Following the optimal sequence, the partial completion corresponds to the best possible subset of comparisons with the same cardinality, independently of when the decision maker stops the filling-in process. The results are confirmed by an experimental study \cite{SzadoczkiBozokiSiposGalambosi2025a}.
The vertex-labelled variant of the graph of graphs approach remains a topic for future research.

\section{Potential areas of applications} \label{Sec7}

In some designs of sports tournaments, players do not play against all the other players due to time or logistical constraints \cite{DevriesereCsatoGoossens2025, DevriesereGoossens2025}. One of them is the Swiss-system, which is widely used in chess.
\cite{Csato2013a} proposes alternative rankings for the 39th Chess Olympiad 2010 Open tournament based on incomplete pairwise comparison matrices. This competition consisted of 11 rounds and was contested by 149 teams. Each match was played on four boards, with three possible results (white win, draw, black win) in each game. This implies nine possible values of game points for each team in every match, which were converted to numerical pairwise comparisons according to four scales. The incomplete logarithmic least squares~\eqref{eq_LLSM} and eigenvector~\eqref{eq_EM} methods were used to compute the weights and rank the teams.

The unique feature of this application is the high number of alternatives (149) and the sparse incomplete pairwise comparison matrix: out of the 11026 elements in the upper triangle, only 810---the total number of matches played by the teams---are known (7.3\%). Nonetheless, solving problem~\eqref{eq_LLSM} remains straightforward as a system of linear equations. On the other hand, the implementation of the incomplete eigenvector method required a small trick. \cite[Section~5]{BozokiFulopRonyai2010} suggests the iterative method of cyclic coordinates to solve problem~\eqref{eq_EM} that can start from any initial vector. While the natural solution chosen by \cite{BozokiFulopRonyai2010} is the all-ones vector, a substantial reduction in computation time can be achieved by starting from the optimal completion corresponding to problem~\eqref{eq_LLSM}, as the weights according to the two weighing procedures are usually quite close.

\cite{Csato2017a} revisits the problem of ranking in Swiss-system chess team tournaments by adopting an axiomatic approach, and analyses the results of two European championships.
However, instead of multiplicative pairwise comparison matrices, the paper uses the additive setting, where the least squares method is equivalent to solving problem~\eqref{eq_LLSM}. This connection is important because \cite{Gonzalez-DiazHendrickxLohmann2014} offers a comprehensive axiomatic study of the least squares method, and \cite{Csato2015a} uncovers how its solution can be obtained as a limit point of an iterative process based on the graph representation of the incomplete pairwise comparison matrix.

While the pairing of the Swiss-system is dynamic \cite{SauerCsehLenzner2024}, incomplete round-robin tournaments have a fixed schedule, that is, the set of pairwise comparisons to be made (matches to be played) is fixed in advance, before the start of the tournament \cite{LiVanBulckGoossens2025}. Ranking in these contests received attention when the outbreak of the Covid-19 pandemic abruptly ended some round-robin tournaments \cite{Csato2021d}.
However, the importance of incomplete round-robin tournaments has greatly increased recently, after UEFA (Union of European Football Associations) decided to organise the most prestigious club football tournaments in the continent, in this format from the 2024/25 season. Ranking in these tournaments, where the set of matches played is given by a regular graph, remains challenging because of the complex interaction between the strength of opponents and the results of the teams \cite{CsatoDevriesereGoossensGyimesiLambersSpieksma2026}.

Incomplete pairwise comparisons emerge in journal and university rankings, too.
\cite{CsatoToth2020} adopt the methodology of \cite{Csato2017a} to rank Hungarian university faculties between 2001 and 2016. Since the centralised assignment system of the Hungarian higher education assigns any student to the first programme where the score limit is reached, the list of applications reveals a partial preference order for each student, from which a sparse incomplete pairwise comparison matrix with a high number of alternatives can be constructed that reflects more than 150 thousand revealed preferences each year.
The resulting rankings represent the collective wisdom of thousands of applicants at a moment of a high-stakes decision. Similar preference-based university rankings are suggested by \cite{AveryGlickmanHoxbyMetrick2013} and \cite{IehleJacqmin2025}.

\cite{WuSuLi2025} uses an analogous idea to evaluate the quality of journals by considering journal lists from several sources: journal $i$ is said to be better than (to win against) journal $j$ if it occupies a higher tier in a given journal list, while they are tied (play a draw) if they are placed in the same category. Naturally, pairwise comparisons may be missing as we do not know whether a journal does not appear in a list due to its low quality or irrelevance to the creator of the list.

Returning to the field of sports, \cite{BozokiCsatoTemesi2016} constructs an incomplete pairwise comparison matrix to evaluate the performance of the 25 leaders in Men's Tennis ATP World Ranking between 1973 and 2013. The pairwise comparison $a_{ij}$ equals $x_{ij} / x_{ji}$, where $x_{ij}$ is the number of matches won by player $i$ against player $j$. The comparison is missing if $x_{ij} + x_{ji} = 0$, that is, players $i$ and $j$ did not play against each other. Finally, if $x_{ij} + x_{ji} > 0$ but $x_{ij} = 0$ ($x_{ji} = 0$), then $a_{ij}$ equals the ceiling $1 / \lceil x_{ji}/5 \rceil$ ($\lceil x_{ij}/5 \rceil$) in the first adjustment, and $1/ \left( x_{ji} + 2 \right)$ ($x_{ij} + 2$) in the second adjustment.

Since the number of matches played by two players can reach 39 (for the pair \emph{Djokovic} vs \emph{Nadal}), another version that accounts for the number of matches played by two players is introduced:
\[
a_{ij} = x_{ij}^{\left( x_{ij} + x_{ji} \right)/39}.
\]
Hence, the numerical comparison of two players remains close to one if they played only a few matches against each other.

The weights are computed by both the incomplete eigenvector and logarithmic least squares methods. While the impact of the adjustment used to avoid division by zero and the procedure for deriving priorities is marginal, the transformation of the pairwise comparison based on the number of matches played is far from negligible. Nonetheless, both the top three (\emph{Nadal}, \emph{Federer}, \emph{Sampras}) and the bottom three players turn out to be robust in the eight rankings.
   
The rankings of 28 top female tennis players \cite{TemesiSzadoczkiBozoki2024} and 14 Go grandmasters together with AlphaGo \cite{ChaoKouLiPeng2018}, a project based on artificial intelligence, follow the same approach by constructing incomplete pairwise comparison matrices based on head-to-head results. The ranking of 1544 Go players demonstrates that large-scale data can be efficiently handled by the incomplete logarithmic least squares method.

A recent innovative ranking of famous violinists \cite{PuppeTasnadi2026} defines the partial rankings of the artists by the number of YouTube views for the same piece played by them. The pairwise comparison of violinists $i$ and $j$ is determined by the number of rankings in which $i$ is ranked above $j$ and $j$ is ranked above $i$, analogous to \cite{BozokiCsatoTemesi2016}. The best performances do not fade over time: listening to violinists from decades ago does not appear to attract less audience than listening to contemporary violinists.

Based on a network meta-analysis of 211 articles, \cite{Voigt2025} ranks 57 removal (customer or route) and 42 insertion operators, used in adaptive large neighborhood search for vehicle routing problems, with respect to their effectiveness. The pairwise comparisons of the operators are derived such that a more frequently used operator wins against a less frequently used one. The incomplete pairwise comparison matrix is constructed analogous to \cite{BozokiCsatoTemesi2016} with the second adjustment, when $a_{ij} = x_{ij} + 2$ if $x_{ji} = 0$. A comparison is missing if two operators are not analysed in the same article.

\section{Conclusions} \label{Sec8}

We have surveyed some fields where incomplete pairwise comparison matrices can be useful, and discussed how rankings can be derived from them. The two most popular weighting procedures, the logarithmic least squares and eigenvector methods, have a natural extension to this general setting (Section~\ref{Sec3}). However, they can give different results, which raises some open questions, such as:
\begin{itemize}
\item
What is the set of incomplete pairwise comparison matrices for which the two approaches give the same estimation of missing entries and/or the same weights?
\item
Does this equivalence hold if the optimal completion is obtained by minimising another inconsistency index?
\item
When does an incomplete pairwise comparison matrix have only one reasonable optimal completion?
\end{itemize}

Section~\ref{Sec4} has illustrated that, in contrast to the complete case, when the logarithmic least squares solution is difficult to debate due to its axiomatic characterisations \cite{Fichtner1984, Barzilai1997, Csato2018b, Csato2019a}, these natural methods can be controversial for incomplete pairwise comparison matrices generated by directed acyclic graphs. This calls for the consideration of novel methods, for instance, the lexicographically optimal completion.

Besides completion and weighting methods, the application of incomplete pairwise comparison matrices may require inconsistency thresholds and a careful choice of questioning order. The order could be especially important in online questionnaires. Even though Sections~\ref{Sec5} and \ref{Sec6} have outlined recent results in both areas, considerable scope remains for future research.
Last but not least, we have seen that there are some common issues in the design of pairwise comparisons and sports tournaments; for example, finding the optimal set of comparisons/matches with a limited budget when each additional comparison or match carries a cost.
In our opinion, both areas of research can benefit from a stronger collaboration.

\begin{acknowledgement}
The research was supported by the National Research, Development and Innovation Office under Grants Advanced 152220 and FK 145838, and the J\'anos Bolyai Research Scholarship of the Hungarian Academy of Sciences.
\end{acknowledgement}

\ethics{Competing Interests}{The authors have no conflicts of interest to declare that are relevant to the content of this chapter.}

\bibliographystyle{spmpsci}
\bibliography{All_references}

\end{document}

%% file: Figure6_optimal_filling_in_sequences.tex
\begin{figure}[ht!]
\centering
\begin{tikzpicture}[every node/.style={circle,inner sep=2pt,draw=black,fill=black!20}, scale=0.69]

\tikzstyle{node2} = [draw=white,fill=white];

%Small nodes
 \node (1A) at (0,0.4) {};
  \node (1B) at (-0.5,-1.6) {};
  \node (1C) at (-1,-0.6) {};
  \node (1D) at (0.5,-1.6) {};
  \node (1E) at (1,-0.6) {};
  
  \node (2A) at (4,0.4) {};
  \node (2B) at (3,-0.6) {};
  \node (2C) at (4.5,-1.6) {};
  \node (2D) at (5,-0.6) {};
  \node (2E) at (3.5,-1.6) {};
  
  \node (3A) at (8.5,-1.6) {};
  \node (3B) at (7.5,-1.6) {};
  \node (3C) at (7,-0.6) {};
  \node (3D) at (8,0.4) {};
  \node (3E) at (9,-0.6) {};
  \node [node2] (4el) at (-5,0.4) {$m=6$};
  
   \node (4A) at (-3,-3.3) {};
  \node (4B) at (-3.5,-5.3) {};
  \node (4C) at (-4,-4.3) {};
  \node (4D) at (-2.5,-5.3) {};
  \node (4E) at (-2,-4.3) {};
  
  \node (5A) at (0.5,-3.3) {};
  \node (5B) at (1.5,-4.3) {};
  \node (5C) at (-0.5,-4.3) {};
  \node (5D) at (1,-5.3) {};
  \node (5E) at (0,-5.3) {};
  
   \node (6A) at (4,-3.3) {};
  \node (6B) at (5,-4.3) {};
  \node (6C) at (4.5,-5.3) {};
  \node (6D) at (3,-4.3) {};
  \node (6E) at (3.5,-5.3) {};
  
  \node (7A) at (6.5,-4.3) {};
  \node (7B) at (7,-5.3) {};
  \node (7C) at (7.5,-3.3) {};
  \node (7D) at (8.5,-4.3) {};
  \node (7E) at (8,-5.3) {};
  
  \node (8A) at (10.5,-3.3) {};
  \node (8B) at (11,-5.3) {};
  \node (8C) at (9.5,-4.3) {};
  \node (8D) at (11.5,-4.3) {};
  \node (8E) at (10,-5.3) {};
  \node [node2] (5el) at (-5,-3.3) {$m=5$};
  
  \node (9A) at (-3,-7.2) {};
  \node (9B) at (-3.5,-9.2) {};
  \node (9C) at (-4,-8.2) {};
  \node (9D) at (-2.5,-9.2) {};
  \node (9E) at (-2,-8.2) {};
  
  \node (10A) at (0.5,-7.2) {};
  \node (10B) at (1,-9.2) {};
  \node (10C) at (-0.5,-8.2) {};
  \node (10D) at (1.5,-8.2) {};
  \node (10E) at (0,-9.2) {};
  
  \node (11A) at (3,-8.2) {};
  \node (11B) at (3.5,-9.2) {};
  \node (11C) at (4,-7.2) {};
  \node (11D) at (5,-8.2) {};
  \node (11E) at (4.5,-9.2) {};
  
  \node (12A) at (8.5,-8.2) {};
  \node (12B) at (7,-9.2) {};
  \node (12C) at (7.5,-7.2) {};
  \node (12D) at (8,-9.2) {};
  \node (12E) at (6.5,-8.2) {};
 
  \node (13A) at (10.5,-7.2) {};
  \node (13B) at (11.5,-8.2) {};
  \node (13C) at (11,-9.2) {};
  \node (13D) at (9.5,-8.2) {};
  \node (13E) at (10,-9.2) {};
  \node [node2] (6el) at (-5,-7.2) {$m=4$};
  
  \node (14A) at (-1,-10.7) {};
  \node (14B) at (-2,-11.7) {};
  \node (14C) at (-0.5,-12.7) {};
  \node (14D) at (0,-11.7) {};
  \node (14E) at (-1.5,-12.7) {};
  
  \node (15A) at (2,-10.7) {};
  \node (15B) at (1.5,-12.7) {};
  \node (15C) at (1,-11.7) {};
  \node (15D) at (2.5,-12.7) {};
  \node (15E) at (3,-11.7) {};

  \node (16A) at (5,-12.7) {};
  \node (16B) at (4.5,-11.7) {};
  \node (16C) at (5.5,-10.7) {};
  \node (16D) at (6,-12.7) {};
  \node (16E) at (6.5,-11.7) {};
  
  \node (17A) at (7.5,-11.7) {};
  \node (17B) at (9,-12.7) {};
  \node (17C) at (8.5,-10.7) {};
  \node (17D) at (8,-12.7) {};
  \node (17E) at (9.5,-11.7) {};
  \node [node2] (7el) at (-5,-10.7) {$m=3$};
  
  \node (18A) at (3,-15.2) {};
  \node (18B) at (2.5,-16.2) {};
  \node (18C) at (1.5,-16.2) {};
  \node (18D) at (2,-14.2) {};
  \node (18E) at (1,-15.2) {};
  
  \node (19A) at (6.5,-15.2) {};
  \node (19B) at (5,-16.2) {};
  \node (19C) at (5.5,-14.2) {};
  \node (19D) at (6,-16.2) {};
  \node (19E) at (4.5,-15.2) {};
  \node [node2] (8el) at (-5,-14.2) {$m=2$};
  
  \node (20A) at (4,-17.2) {};
  \node (20B) at (3.5,-19.2) {};
  \node (20C) at (3,-18.2) {};
  \node (20D) at (5,-18.2) {};
  \node (20E) at (4.5,-19.2) {};
  \node [node2] (9el) at (-5,-17.2) {$m=1$};
  
  \node (21A) at (4,-20.2) {};
  \node (21B) at (4.5,-22.2) {};
  \node (21C) at (3,-21.2) {};
  \node (21D) at (5,-21.2) {};
  \node (21E) at (3.5,-22.2) {};
  \node [node2] (10el) at (-5,-20.2) {$m=0$};

%Small edges
  \draw (1A) -- (1B)
        (1A) -- (1C)
        (1A) -- (1D)
        (1A) -- (1E)
        
        (2A) -- (2B)
        (2A) -- (2C)
        (2A) -- (2D)
        (2C) -- (2E)
        
        (3A) -- (3B)
        (3B) -- (3C)
        (3C) -- (3D)
        (3D) -- (3E)
        
        (4A) -- (4B)
        (4A) -- (4C)
        (4A) -- (4D)
        (4A) -- (4E)
        (4B) -- (4D)
        
        (5A) -- (5B)
        (5A) -- (5C)
        (5A) -- (5D)
        (5C) -- (5E)
        (5D) -- (5E)
        
        (6A) -- (6B)
        (6A) -- (6C)
        (6A) -- (6D)
        (6C) -- (6E)
        (6C) -- (6D)
        
        (7B) -- (7E)
        (7A) -- (7C)
        (7A) -- (7D)
        (7C) -- (7E)
        (7C) -- (7D)
        
        (8B) -- (8E)
        (8A) -- (8C)
        (8A) -- (8D)
        (8C) -- (8E)
        (8B) -- (8D)
        
        (9A) -- (9B)
        (9A) -- (9C)
        (9A) -- (9D)
        (9A) -- (9E)
        (9B) -- (9D)
        (9B) -- (9C)
        
        (10A) -- (10B)
        (10A) -- (10C)
        (10A) -- (10D)
        (10C) -- (10E)
        (10D) -- (10E)
        (10B) -- (10E)
        
        (11B) -- (11E)
        (11A) -- (11C)
        (11A) -- (11D)
        (11C) -- (11E)
        (11C) -- (11D)
        (11C) -- (11B)
        
        (12B) -- (12E)
        (12A) -- (12C)
        (12A) -- (12D)
        (12C) -- (12E)
        (12C) -- (12D)
        (12B) -- (12D)
        
        (13A) -- (13B)
        (13A) -- (13C)
        (13A) -- (13D)
        (13C) -- (13E)
        (13C) -- (13D)
        (13B) -- (13D)
        
        (14A) -- (14B)
        (14A) -- (14C)
        (14A) -- (14D)
        (14C) -- (14E)
        (14D) -- (14E)
        (14B) -- (14E)
        (14A) -- (14E)
        
        (15A) -- (15E)
        (15A) -- (15C)
        (15A) -- (15D)
        (15D) -- (15E)
        (15B) -- (15D)
        (15A) -- (15B)
        (15C) -- (15B)
                
        (16A) -- (16B)
        (16A) -- (16C)
        (16A) -- (16D)
        (16C) -- (16E)
        (16C) -- (16D)
        (16B) -- (16D)
        (16C) -- (16B)
        
        (17B) -- (17E)
        (17A) -- (17C)
        (17A) -- (17D)
        (17C) -- (17E)
        (17C) -- (17B)
        (17B) -- (17D)
        (17D) -- (17E)
        
        (18B) -- (18E)
        (18A) -- (18C)
        (18A) -- (18D)
        (18C) -- (18E)
        (18C) -- (18B)
        (18B) -- (18D)
        (18D) -- (18E)
        (18C) -- (18D)
        
        (19B) -- (19E)
        (19A) -- (19C)
        (19A) -- (19D)
        (19C) -- (19E)
        (19B) -- (19A)
        (19B) -- (19D)
        (19A) -- (19E)
        (19C) -- (19D)
        
        (20B) -- (20E)
        (20A) -- (20C)
        (20A) -- (20D)
        (20C) -- (20E)
        (20C) -- (20B)
        (20B) -- (20D)
        (20D) -- (20E)
        (20C) -- (20D)
        (20A) -- (20E)
        
        (21B) -- (21E)
        (21A) -- (21C)
        (21A) -- (21D)
        (21C) -- (21E)
        (21B) -- (21A)
        (21B) -- (21D)
        (21A) -- (21E)
        (21C) -- (21D)
        (21D) -- (21E)
        (21B) -- (21C);

%Large edges   
    \draw[very thick] (0,-1.75) -- (-3,-3.1)
                      (4,-1.7) -- (-3,-3.1)
                      (4,-1.7) -- (0.5,-3.1)
                      (4,-1.7) -- (4,-3.1)
                      (4,-1.7) -- (7.5,-3.1)
                      (8,-1.7) -- (0.5,-3.1)
                      (8,-1.7) -- (4,-3.1)
                      (8,-1.7) -- (7.5,-3.1)
                      (8,-1.7) -- (10.5,-3.1)
                      (-3,-5.5) -- (-3,-7)
                      (-3,-5.5) -- (4,-7)
                      (0.5,-5.5) -- (-3,-7)
                      (0.5,-5.5) -- (0.5,-7)
                      (0.5,-5.5) -- (7.5,-7)
                      (0.5,-5.5) -- (10.5,-7)
                      (4,-5.5) -- (-3,-7)
                      (4,-5.5) -- (7.5,-7)
                      (4,-5.5) -- (10.5,-7)
                      (7.5,-5.5) -- (4,-7)
                      (7.5,-5.5) -- (7.5,-7)
                      (7.5,-5.5) -- (10.5,-7)
                      (10.5,-5.5) -- (7.5,-7)
                      (-3,-9.4) -- (-1,-10.4)
                      (-3,-9.4) -- (2,-10.4)
                      (-3,-9.4) -- (5.5,-10.4)
                      (0.5,-9.35) -- (-1,-10.4)
                      (4,-9.4) -- (2,-10.4)
                      (7.5,-9.4) -- (2,-10.4)
                      (7.5,-9.4) -- (8.5,-10.4)
                      (10.5,-9.4) -- (2,-10.4)
                      (10.5,-9.4) -- (5.5,-10.4)
                      (10.5,-9.4) -- (8.5,-10.4)
                      (-1,-12.9) -- (2.5,-14)
                      (2,-12.9) -- (2.5,-14)
                      (5.5,-12.9) -- (2.5,-14)
                      (8.5,-12.9) -- (2.5,-14)
                      (2,-12.9) -- (6,-14)
                      (2,-16.4) -- (4,-17);
                      
%Rectangles
\draw[line width=0.75mm, draw=green,fill=green!50,opacity=0.3] (-1.2,-1.75) rectangle (1.2,0.9);
\draw[very thick, draw=black!50,fill=black!50,opacity=0.08] (2.8,-1.7) rectangle (5.2,0.8);
\draw[very thick, draw=black!50,fill=black!50,opacity=0.08] (6.8,-1.7) rectangle (9.2,0.9);
\draw[very thick, draw=black!50,fill=black!50,opacity=0.08] (-4.2,-5.5) rectangle (-1.8,-3.1);
\draw[very thick, draw=black!50,fill=black!50,opacity=0.08] (-0.7,-5.5) rectangle (1.7,-3.1);
\draw[very thick, draw=black!50,fill=black!50,opacity=0.08] (2.8,-5.5) rectangle (5.2,-3.1);
\draw[very thick, draw=black!50,fill=black!50,opacity=0.08] (6.3,-5.5) rectangle (8.7,-3.1);
\draw[line width=0.75mm, draw=green,fill=green!50,opacity=0.3] (9.3,-5.5) rectangle (11.7,-3.1);
\draw[very thick, draw=black!50,fill=black!50,opacity=0.08] (-4.2,-9.4) rectangle (-1.8,-7);
\draw[line width=0.75mm, draw=green,fill=green!50,opacity=0.3] (-0.7,-9.35) rectangle (1.7,-7);
\draw[very thick, draw=black!50,fill=black!50,opacity=0.08] (2.8,-9.4) rectangle (5.2,-7);
\draw[very thick, draw=black!50,fill=black!50,opacity=0.08] (6.3,-9.4) rectangle (8.7,-7);
\draw[very thick, draw=black!50,fill=black!50,opacity=0.08] (9.3,-9.4) rectangle (11.7,-7);
\draw[very thick, draw=black!50,fill=black!50,opacity=0.08] (-2.2,-12.9) rectangle (0.2,-10.4);
\draw[very thick, draw=black!50,fill=black!50,opacity=0.08] (0.8,-12.9) rectangle (3.2,-10.4);
\draw[very thick, draw=black!50,fill=black!50,opacity=0.08] (4.3,-12.9) rectangle (6.7,-10.4);
\draw[line width=0.75mm, draw=green,fill=green!50,opacity=0.3] (7.3,-12.9) rectangle (9.7,-10.4);
\draw[very thick, draw=black!50,fill=black!50,opacity=0.08] (0.8,-16.4) rectangle (3.2,-14);
\draw[line width=0.75mm, draw=green,fill=green!50,opacity=0.3] (4.2,-16.4) rectangle (6.7,-14);
\draw[line width=0.75mm, draw=green,fill=green!50,opacity=0.3] (2.8,-19.4) rectangle (5.2,-17);
\draw[line width=0.75mm, draw=green,fill=green!50,opacity=0.3] (2.7,-22.4) rectangle (5.2,-20);

\draw[very thick, draw=green]
                      (8.5,-12.9) -- (6,-14)
                      (5.5,-16.4) -- (4,-17)
                      (4,-19.4) -- (4,-20)
                      (0.5,-9.35) -- (8.5,-10.4);
\end{tikzpicture}
\caption{Optimal filling-in patterns and sequence ($n=5$, optimal graphs are highlighted by green)}
\label{Fig6}
\end{figure}
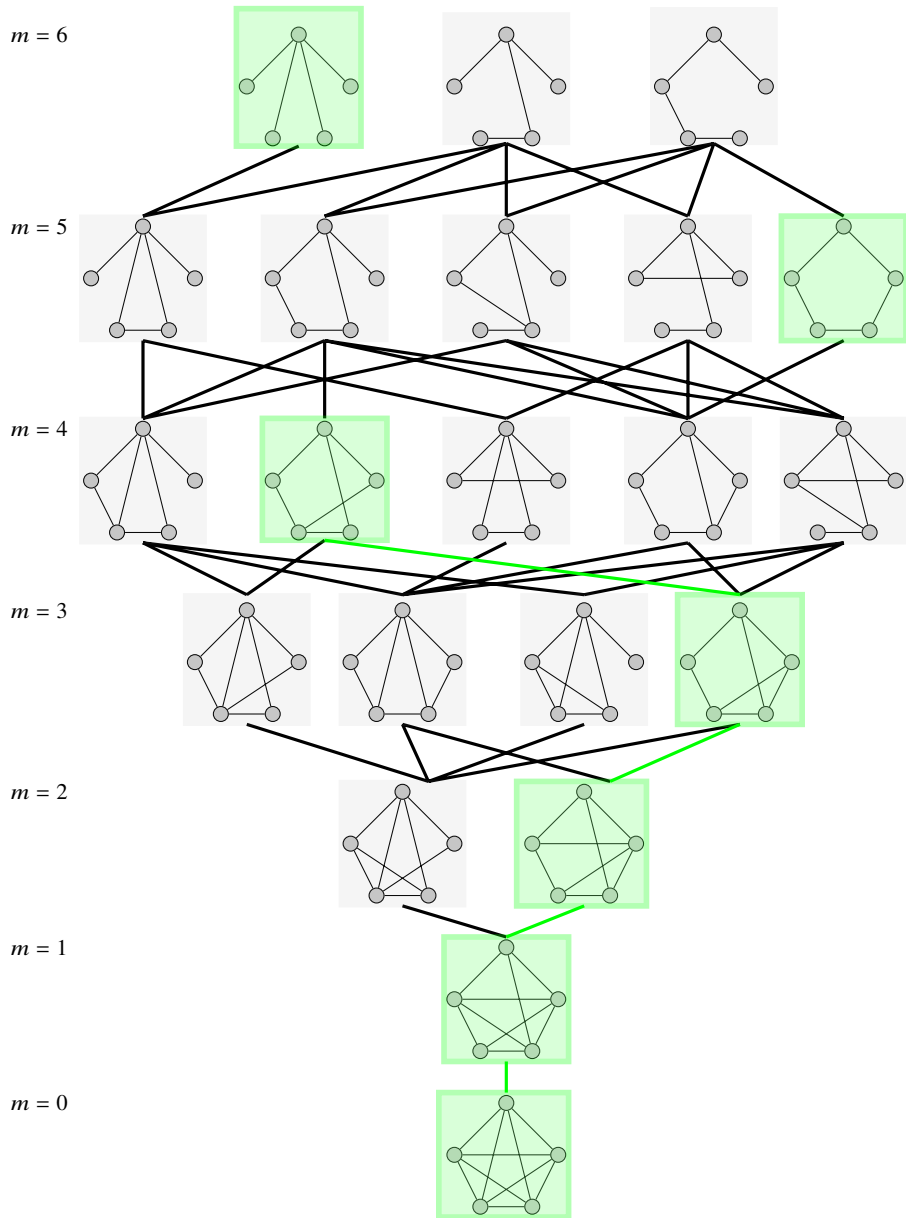

%\end{document}

%% file: Incomplete_PCM_theory_applications.bbl
\begin{thebibliography}{10}
\providecommand{\url}[1]{{#1}}
\providecommand{\urlprefix}{URL }
\expandafter\ifx\csname urlstyle\endcsname\relax
  \providecommand{\doi}[1]{DOI~\discretionary{}{}{}#1}\else
  \providecommand{\doi}{DOI~\discretionary{}{}{}\begingroup
  \urlstyle{rm}\Url}\fi

\bibitem{AbastanteCorrenteGrecoIsizakaLami2019}
Abastante, F., Corrente, S., Greco, S., Ishizaka, A., Lami, I.M.: A new
  parsimonious {AHP} methodology: Assigning priorities to many objects by
  comparing pairwise few reference objects.
\newblock Expert Systems with Applications \textbf{127}, 109--120 (2019)

\bibitem{AgostonCsato2022}
{\'A}goston, K.{\relax Cs}., Csat{\'o}, L.: Inconsistency thresholds for
  incomplete pairwise comparison matrices.
\newblock Omega \textbf{108}, 102576 (2022)

\bibitem{AgostonCsato2024}
{\'A}goston, K.{\relax Cs}., Csat{\'o}, L.: A lexicographically optimal
  completion for pairwise comparison matrices with missing entries.
\newblock European Journal of Operational Research \textbf{314}(3), 1078--1086
  (2024)

\bibitem{AgostonCsato2026}
{\'A}goston, K.{\relax Cs}., Csat{\'o}, L.: Refined thresholds for
  inconsistency: The effect of the graph associated with incomplete pairwise
  comparisons.
\newblock Expert Systems with Applications \textbf{328}, 132938 (2026)

\bibitem{AguaronMoreno-Jimenez2003}
Aguar\'on, J., Moreno-Jim{\'e}nez, J.M.: The geometric consistency index:
  Approximated thresholds.
\newblock European Journal of Operational Research \textbf{147}(1), 137--145
  (2003)

\bibitem{AlonsoLamata2006}
Alonso, J.A., Lamata, M.T.: Consistency in the {A}nalytic {H}ierarchy
  {P}rocess: {A} new approach.
\newblock International Journal of Uncertainty, Fuzziness and Knowledge-Based
  Systems \textbf{14}(4), 445--459 (2006)

\bibitem{AveryGlickmanHoxbyMetrick2013}
Avery, C.N., Glickman, M.E., Hoxby, C.M., Metrick, A.: A revealed preference
  ranking of {U}.{S}.~colleges and universities.
\newblock Quarterly Journal of Economics \textbf{128}(1), 425--467 (2013)

\bibitem{Barzilai1997}
Barzilai, J.: Deriving weights from pairwise comparison matrices.
\newblock Journal of the Operational Research Society \textbf{48}(12),
  1226--1232 (1997)

\bibitem{BozokiCsatoTemesi2016}
Boz{\'o}ki, S., Csat{\'o}, L., Temesi, J.: An application of incomplete
  pairwise comparison matrices for ranking top tennis players.
\newblock European Journal of Operational Research \textbf{248}(1), 211--218
  (2016)

\bibitem{BozokiDezsoPoeszTemesi2013}
Boz{\'o}ki, S., Dezs{\H{o}}, L., Poesz, A., Temesi, J.: Analysis of pairwise
  comparison matrices: an empirical research.
\newblock Annals of Operations Research \textbf{211}(1), 511--528 (2013)

\bibitem{BozokiFulopRonyai2010}
Boz\'oki, S., F\"ul\"op, J., R\'onyai, L.: On optimal completion of incomplete
  pairwise comparison matrices.
\newblock Mathematical and Computer Modelling \textbf{52}(1-2), 318--333 (2010)

\bibitem{BozokiTsyganok2019}
Boz{\'o}ki, S., Tsyganok, V.: The (logarithmic) least squares optimality of the
  arithmetic (geometric) mean of weight vectors calculated from all spanning
  trees for incomplete additive (multiplicative) pairwise comparison matrices.
\newblock International Journal of General Systems \textbf{48}(4), 362--381
  (2019)

\bibitem{Brunelli2018}
Brunelli, M.: A survey of inconsistency indices for pairwise comparisons.
\newblock International Journal of General Systems \textbf{47}(8), 751--771
  (2018)

\bibitem{CarmoneKaraZanakis1997}
Carmone, F.J.J., Kara, A., Zanakis, S.H.: A {M}onte {C}arlo investigation of
  incomplete pairwise comparison matrices in {AHP}.
\newblock European Journal of Operational Research \textbf{102}(3), 538--553
  (1997)

\bibitem{Cavallo2020}
Cavallo, B.: Functional relations and {S}pearman correlation between
  consistency indices.
\newblock Journal of the Operational Research Society \textbf{71}(2), 301--311
  (2020)

\bibitem{ChaoKouLiPeng2018}
Chao, X., Kou, G., Li, T., Peng, Y.: Jie {K}e versus {A}lpha{G}o: A ranking
  approach using decision making method for large-scale data with incomplete
  information.
\newblock European Journal of Operational Research \textbf{265}(1), 239--247
  (2018)

\bibitem{CrawfordWilliams1985}
Crawford, G., Williams, C.: A note on the analysis of subjective judgment
  matrices.
\newblock Journal of Mathematical Psychology \textbf{29}(4), 387--405 (1985)

\bibitem{Csato2013a}
Csat{\'o}, L.: Ranking by pairwise comparisons for {S}wiss-system tournaments.
\newblock Central European Journal of Operations Research \textbf{21}(4),
  783--803 (2013)

\bibitem{Csato2015a}
Csat{\'o}, L.: A graph interpretation of the least squares ranking method.
\newblock Social Choice and Welfare \textbf{44}(1), 51--69 (2015)

\bibitem{Csato2017a}
Csat{\'o}, L.: On the ranking of a {S}wiss system chess team tournament.
\newblock Annals of Operations Research \textbf{254}(1-2), 17--36 (2017)

\bibitem{Csato2018a}
Csat{\'o}, L.: Characterization of an inconsistency ranking for pairwise
  comparison matrices.
\newblock Annals of Operations Research \textbf{261}(1-2), 155--165 (2018)

\bibitem{Csato2018b}
Csat{\'o}, L.: Characterization of the row geometric mean ranking with a group
  consensus axiom.
\newblock Group Decision and Negotiation \textbf{27}(6), 1011--1027 (2018)

\bibitem{Csato2019c}
Csat{\'o}, L.: Axiomatizations of inconsistency indices for triads.
\newblock Annals of Operations Research \textbf{280}(1-2), 99--110 (2019)

\bibitem{Csato2019a}
Csat{\'o}, L.: A characterization of the {L}ogarithmic {L}east {S}quares
  {M}ethod.
\newblock European Journal of Operational Research \textbf{276}(1), 212--216
  (2019)

\bibitem{Csato2021d}
Csat{\'o}, L.: Coronavirus and sports leagues: obtaining a fair ranking when
  the season cannot resume.
\newblock IMA Journal of Management Mathematics \textbf{32}(4), 547--560 (2021)

\bibitem{Csato2024b}
Csat{\'o}, L.: How to choose a completion method for pairwise comparison
  matrices with missing entries: {A}n axiomatic result.
\newblock International Journal of Approximate Reasoning \textbf{164}, 109063
  (2024)

\bibitem{Csato2025a}
Csat{\'o}, L.: The logarithmic least squares priorities and ordinal violations
  in the best-worst method.
\newblock Expert Systems with Applications \textbf{265}, 125966 (2025)

\bibitem{CsatoAgostonBozoki2024}
Csat{\'o}, L., {\'A}goston, K.{\relax Cs}., Boz{\'o}ki, S.: On the coincidence
  of optimal completions for small pairwise comparison matrices with missing
  entries.
\newblock Annals of Operations Research \textbf{331}(1), 239--247 (2024)

\bibitem{CsatoDevriesereGoossensGyimesiLambersSpieksma2026}
Csat{\'o}, L., Devriesere, K., Goossens, D., Gyimesi, A., Lambers, R.,
  Spieksma, F.: Ranking matters: Does the new format select the best teams for
  the knockout phase in the {UEFA} {C}hampions {L}eague?
\newblock International Journal of Sports Science \& Coaching \textbf{21}(2),
  1123–1131 (2026)

\bibitem{CsatoRonyai2016}
Csat{\'o}, L., R{\'o}nyai, L.: Incomplete pairwise comparison matrices and
  weighting methods.
\newblock Fundamenta Informaticae \textbf{144}(3-4), 309--320 (2016)

\bibitem{CsatoToth2020}
Csat{\'o}, L., T{\'o}th, {\relax Cs}.: University rankings from the revealed
  preferences of the applicants.
\newblock European Journal of Operational Research \textbf{286}(1), 309--320
  (2020)

\bibitem{DeGraan1980}
De~Graan, J.G.: Extensions of the multiple criteria analysis method of {T}.
  {L}. {S}aaty.
\newblock Report, National Institute for Water Supply, Voorburg (1980)

\bibitem{DevriesereCsatoGoossens2025}
Devriesere, K., Csat{\'o}, L., Goossens, D.: Tournament design: A review from
  an operational research perspective.
\newblock European Journal of Operational Research \textbf{324}(1), 1--21
  (2025)

\bibitem{DevriesereGoossens2025}
Devriesere, K., Goossens, D.: Redesigning {B}elgian youth field hockey
  competitions using an incomplete round-robin tournament.
\newblock INFORMS Journal on Applied Analytics \textbf{55}(6), 457--468 (2025)

\bibitem{DuszakKoczkodaj1994}
Duszak, Z., Koczkodaj, W.W.: Generalization of a new definition of consistency
  for pairwise comparisons.
\newblock Information Processing Letters \textbf{52}(5), 273--276 (1994)

\bibitem{Fichtner1984}
Fichtner, J.: Some thoughts about the mathematics of the {A}nalytic {H}ierarchy
  {P}rocess.
\newblock Tech. rep., Institut f{\"u}r Angewandte Systemforschung und
  Operations Research, Universit\"at der Bundeswehr M\"unchen (1984)

\bibitem{Gonzalez-DiazHendrickxLohmann2014}
Gonz\'alez-D\'iaz, J., Hendrickx, R., Lohmann, E.: Paired comparisons analysis:
  an axiomatic approach to ranking methods.
\newblock Social Choice and Welfare \textbf{42}(1), 139--169 (2014)

\bibitem{Harker1987a}
Harker, P.T.: Alternative modes of questioning in the {A}nalytic {H}ierarchy
  {P}rocess.
\newblock Mathematical Modelling \textbf{9}(3-5), 353--360 (1987)

\bibitem{IehleJacqmin2025}
Iehl{\'e}, V., Jacqmin, J.: Can students help us better rank higher education
  institutions? {T}he case of {F}rench business schools.
\newblock Studies in Higher Education \textbf{50}(12), 2857--2880 (2025)

\bibitem{deJong1984}
de~Jong, P.: A statistical approach to {S}aaty's scaling method for priorities.
\newblock Journal of Mathematical Psychology \textbf{28}(4), 467--478 (1984)

\bibitem{KaiserSerlin1978}
Kaiser, H.F., Serlin, R.C.: Contributions to the method of paired comparisons.
\newblock Applied Psychological Measurement \textbf{2}(3), 423--432 (1978)

\bibitem{Koczkodaj1993}
Koczkodaj, W.W.: A new definition of consistency of pairwise comparisons.
\newblock Mathematical and Computer Modelling \textbf{18}(7), 79--84 (1993)

\bibitem{Kulakowski2015b}
Ku{\l}akowski, K.: A heuristic rating estimation algorithm for the pairwise
  comparisons method.
\newblock Central European Journal of Operations Research \textbf{23},
  187–203 (2015)

\bibitem{KulakowskiTalaga2020}
Ku{\l}akowski, K., Talaga, D.: Inconsistency indices for incomplete pairwise
  comparisons matrices.
\newblock International Journal of General Systems \textbf{49}(2), 174--200
  (2020)

\bibitem{Kwiesielewicz1996}
Kwiesielewicz, M.: The logarithmic least squares and the generalized
  pseudoinverse in estimating ratios.
\newblock European Journal of Operational Research \textbf{93}(3), 611--619
  (1996)

\bibitem{LiVanBulckGoossens2025}
Li, M., Van~Bulck, D., Goossens, D.: Beyond leagues: A single incomplete round
  robin tournament for multi-league sports timetabling.
\newblock European Journal of Operational Research \textbf{323}(1), 208--223
  (2025)

\bibitem{LundySirajGreco2017}
Lundy, M., Siraj, S., Greco, S.: The mathematical equivalence of the ``spanning
  tree'' and row geometric mean preference vectors and its implications for
  preference analysis.
\newblock European Journal of Operational Research \textbf{257}(1), 197--208
  (2017)

\bibitem{PuppeTasnadi2026}
Puppe, C., Tasn{\'a}di, A.: Do performances of violin virtuosi fade over time?
  {A} response using the {N}ash collective utility function.
\newblock European Journal of Operational Research \textbf{333}(1), 203--215
  (2026)

\bibitem{Rabinowitz1976}
Rabinowitz, G.: Some comments on measuring world influence.
\newblock Conflict Management and Peace Science \textbf{2}(1), 49--55 (1976)

\bibitem{Rezaei2015}
Rezaei, J.: Best-worst multi-criteria decision-making method.
\newblock Omega \textbf{53}, 49--57 (2015)

\bibitem{Ross1934}
Ross, R.T.: Optimum orders for the presentation of pairs in the method of
  paired comparisons.
\newblock Journal of Educational Psychology \textbf{25}(5), 375--382 (1934)

\bibitem{Saaty1977}
Saaty, T.L.: A scaling method for priorities in hierarchical structures.
\newblock Journal of Mathematical Psychology \textbf{15}(3), 234--281 (1977)

\bibitem{Saaty1980}
Saaty, T.L.: The {A}nalytic {H}ierarchy {P}rocess: Planning, Priority Setting,
  Resource Allocation.
\newblock McGraw-Hill, New York, New York, USA (1980)

\bibitem{SaatyOzdemir2003}
Saaty, T.L., Ozdemir, M.S.: Why the magic number seven plus or minus two.
\newblock Mathematical and Computer Modelling \textbf{38}(3-4), 233--244 (2003)

\bibitem{SauerCsehLenzner2024}
Sauer, P., Cseh, {\'A}., Lenzner, P.: Improving ranking quality and fairness in
  {S}wiss-system chess tournaments.
\newblock Journal of Quantitative Analysis in Sports \textbf{20}(2), 127--146
  (2024)

\bibitem{Schmeidler1969}
Schmeidler, D.: The nucleolus of a characteristic function game.
\newblock SIAM Journal on Applied Mathematics \textbf{17}(6), 1163--1170 (1969)

\bibitem{ShiraishiObata2002}
Shiraishi, S., Obata, T.: On a maximization problem arising from a positive
  reciprocal matrix in {AHP}.
\newblock Bulletin of Informatics and Cybernetics \textbf{34}(2), 91--96 (2002)

\bibitem{ShiraishiObataDaigo1998}
Shiraishi, S., Obata, T., Daigo, M.: Properties of a positive reciprocal matrix
  and their application to {AHP}.
\newblock Journal of the Operations Research Society of Japan \textbf{41}(3),
  404--414 (1998)

\bibitem{SzadoczkiBozoki2025}
Sz{\'a}doczki, {\relax Zs}., Boz{\'o}ki, S.: Optimal sequences for pairwise
  comparisons: the graph of graphs approach.
\newblock Annals of Operations Research \textbf{353}(3), 1099–1122 (2025)

\bibitem{SzadoczkiBozokiSiposGalambosi2025a}
Sz{\'a}doczki, {\relax Zs}., Boz{\'o}ki, S., Sipos, L., Galambosi, {\relax
  Zs}.: An experimental approach: Converting verbal expressions to numerical
  scales (2025).
\newblock Manuscript. {DOI}:
  \href{https://arxiv.org/abs/2507.04539}{10.48550/arXiv.2507.04539}

\bibitem{SzadoczkiBozokiTekile2022}
Sz{\'a}doczki, {\relax Zs}., Boz{\'o}ki, S., Tekile, H.A.: Filling in pattern
  designs for incomplete pairwise comparison matrices: ({Q}uasi-)regular graphs
  with minimal diameter.
\newblock Omega \textbf{107}, 102557 (2022)

\bibitem{TakedaYu1995}
Takeda, E., Yu, P.L.: Assessing priority weights from subsets of pairwise
  comparisons in multiple criteria optimization problems.
\newblock European Journal of Operational Research \textbf{86}(2), 315--331
  (1995)

\bibitem{TekileBrunelliFedrizzi2023}
Tekile, H.A., Brunelli, M., Fedrizzi, M.: A numerical comparative study of
  completion methods for pairwise comparison matrices.
\newblock Operations Research Perspectives \textbf{10}, 100272 (2023)

\bibitem{TemesiSzadoczkiBozoki2024}
Temesi, J., Sz{\'a}doczki, {\relax Zs}., Boz{\'o}ki, S.: Incomplete pairwise
  comparison matrices: Ranking top women tennis players.
\newblock Journal of the Operational Research Society \textbf{75}(1), 145--157
  (2024)

\bibitem{Tsyganok2000}
Tsyganok, V.: Combinatorial method of pairwise comparisons with feedback.
\newblock Data Recording, Storage \& Processing \textbf{2}, 92--102 (2000)

\bibitem{Tsyganok2010}
Tsyganok, V.: Investigation of the aggregation effectiveness of expert
  estimates obtained by the pairwise comparison method.
\newblock Mathematical and Computer Modelling \textbf{52}(3-4), 538--544 (2010)

\bibitem{TuWuPedrycz2023}
Tu, J., Wu, Z., Pedrycz, W.: Priority ranking for the best-worst method.
\newblock Information Sciences \textbf{635}, 42--55 (2023)

\bibitem{Voigt2025}
Voigt, S.: A review and ranking of operators in adaptive large neighborhood
  search for vehicle routing problems.
\newblock European Journal of Operational Research \textbf{322}(2), 357--375
  (2025)

\bibitem{WilliamsCrawford1980}
Williams, C., Crawford, G.: Analysis of subjective judgment matrices.
\newblock Interim report R-2572-AF, Rand Corporation, Santa Monica (1980)

\bibitem{WuSuLi2025}
Wu, D., Su, Q., Li, J.: Identification of home bias in journal ranking lists.
\newblock Journal of Informetrics \textbf{19}(3), 101707 (2025)

\bibitem{XuWang2024}
Xu, Y., Wang, D.: Some methods to derive the priority weights from the
  best--worst method matrix and weight efficiency test in view of incomplete
  pairwise comparison matrix.
\newblock Fuzzy Optimization and Decision Making \textbf{23}(1), 31--62 (2024)

\end{thebibliography}
